\documentclass[english]{elsarticle}
\usepackage[T1]{fontenc}
\usepackage[latin9]{inputenc}
\usepackage{listings}
\usepackage{babel}
\usepackage{amsmath}
\usepackage{graphicx}
\usepackage[unicode=true]
 {hyperref}

\makeatletter

\providecommand{\tabularnewline}{\\}

\usepackage{url}

\usepackage{listings}
\usepackage{color}

\definecolor{mygreen}{rgb}{0,0.6,0}
\definecolor{mygray}{rgb}{0.5,0.5,0.5}
\definecolor{mymauve}{rgb}{0.58,0,0.82}

\makeatother

\begin{document}
\title{A modern resistive magnetohydrodynamics solver using C++ and the Boost library.\tnoteref{label1}} \tnotetext[label1]{This work was  supported by the Marshall plan scholarship of the Austrian Marshall plan foundation (http://www.marshallplan.at/) and by the Austrian Science Fund (FWF) -- project id: P25346. Computational resources have been provided by the Austrian Ministry of Science BMWF as part of the UniInfrastrukturprogramm of the Focal Point Scientific Computing at the University of Innsbruck. The computational results presented have been achieved [in part] using the Vienna Scientific Cluster (VSC).}
\author[uibk]{Lukas Einkemmer
} \ead{lukas.einkemmer@uibk.ac.at}
\address[uibk]{Department of Mathematics, University of Innsbruck, Austria}
\begin{abstract}
In this paper we describe the implementation of our C++ resistive magnetohydrodynamics solver. The framework developed facilitates the separation of the code implementing the specific numerical method and the physical model, on the one hand, from the handling of boundary conditions and the management of the computational domain, on the other hand. In particular, this will allow us to use finite difference stencils which are only defined in the interior of the domain (the boundary conditions are handled automatically). We will discuss this and other design considerations and their impact on performance in some detail. In addition, we provide a documentation of the code developed and demonstrate that a  performance comparable to Fortran can be achieved, while still maintaining a maximum of code readability and extensibility.
\end{abstract} 
\maketitle

\section{Introduction}

The use of numerical simulations to study problems from the sciences
is prevalent in both industry and academia. Due to the continuous
increase in computing power more realistic (and thus more complex)
physical models are considered in such simulations. This often requires
the adaption of both the physical models as well as the numerical
methods used to solve a given problem.

However, contrary to application development, many scientific codes
are still written either in the Fortran or the C programming language
(although there is a increasing tendency to use selected features
from C++; see, for example, \citep{pesch2007}). In our opinion, the
reason for this is at least threefold. First, performance considerations
are one of the (if not the) primary concern when developing scientific
codes. This explains why interpreted languages (such as Python or
Matlab) have generally only been employed for very specific applications.
Second, object oriented programming has not been as successful as
in other domains (even though, in principle, the value of this abstraction
has been recognized; see, for example, \citep{cary1997}). Third,
the need to interface with legacy code (which often has been extensively
optimized) is an important consideration. Such code is almost exclusively
written in either Fortran or C (many libraries provide interfaces
for both languages). 

Nevertheless, due to increasing code complexity, the need for suitable
mechanisms of abstraction is a major concern if maintainable scientific
codes are to be developed. However, due to performance considerations
many programming paradigms present in the C++ language have either
been dismissed or only used cautiously (such as abstract base classes)
while other where often considered to be too cumbersome to use (such
as templates). It has long been argued that Fortran gives superior
performance to C since the impossibility of pointer aliasing provides
the compiler with more opportunity to generate efficient code. This
argument has to be considered in the context that the C++ standard
library does not include a multidimensional array type. However, template
based libraries have been developed that provide C++ with multidimensional
arrays that rival the performance of Fortran (see e.g. \citep{veldhuizen1998}).

In this paper we will introduce a numerical solver for the equations
of resistive magnetohydrodynamics (resistive MHD). This code employs
a finite difference or finite volume discretization in space (basically
the only assumption with respect to the space discretization made
is that only a single value is stored per cell). For time integration
the CVODE library (which is based on the backward differentiation
formulas) is employed. We will explain the fundamental design principles
of the code in section \ref{sec:Design-principles}. These allow us
to implement a given space or time discretization in a few lines without
considering the specific boundary conditions and to implement the
physical model in only a few lines of code that resembles the mathematical
form of the equations. Nevertheless, binary compatibility with respect
to the CVODE vector type (or any other C or Fortran vector type) is
maintained.

This will be accomplished by the limited use of templates and by using
the Boost library (where appropriate). Our general philosophy is that
we limit ourselves to features of the C++ language for which there
is no clear evidence that a significant performance penalty is incurred.
In addition, we choose a more commonly employed class over a less
well known one if the less well known implementation does not provide
a clear performance advantage. 

Let us also note that the framework we develop contains some template
code. However, this code should be discernible for scientists even
if they only have a passing familiarity with template techniques (in
contrast to libraries such as Boost \citep{boost} or Blitz++ where
significantly more involved template techniques are employed in order
to provide a fast implementation as well as a convenient interface
to the user). Furthermore, to implement a different numerical scheme
or physical model no knowledge of templates techniques is necessary. 

In section \ref{sec:Numerical-examples} we will consider two plasma
physics problems. Numerical simulations and weak scaling tests will
be conducted for the reconnection problem and the Kelvin--Helmholtz
instability. In addition, we also present the results of a strong
scaling test. In section \ref{sec:Documentation} we document the
structure of the source code. Finally, we conclude in section \ref{sec:Conclusion-&-Outlook}.

It is our hope that the framework developed here can be used as a
reference implementation and as a stepping stone if the effect of
different numerical methods or physical models need to be evaluated.

\section{Magnetohydrodynamics \& numerical considerations}

In this section we describe the equations of magnetohydrodynamics
(MHD) which are often used to model a plasma system. These equations
are justified if the assumption is made that the plasma under consideration
is in thermodynamic equilibrium. The equations of MHD are described
in terms of the density $\rho$, the fluid velocity $\boldsymbol{v}$,
the magnetic field $\boldsymbol{B}$, the energy $e$, and the pressure
$p$. 

For the purpose of performing the spatial discretization, these equations
are often cast into the so-called divergence form. Then, the equations
of motion read as (see e.g. \citep{reynolds2010})
\begin{align*}
\frac{\partial U}{\partial t}+\nabla\cdot F(U) & =0
\end{align*}
with state vector
\[
U=\left[\begin{array}{c}
\rho\\
\boldsymbol{v}\\
\boldsymbol{B}\\
e
\end{array}\right]
\]
and
\begin{equation}
F(U)=\left[\begin{array}{c}
\rho\boldsymbol{v}\\
\rho\boldsymbol{v}\otimes\boldsymbol{v}+(p+\tfrac{1}{2}B^{2})I-\boldsymbol{B}\otimes\boldsymbol{B}\\
\boldsymbol{v}\otimes\boldsymbol{B}-\boldsymbol{B}\otimes\boldsymbol{v}\\
(e+p+\tfrac{1}{2}B^{2})\boldsymbol{v}-\boldsymbol{B}(\boldsymbol{B}\cdot\boldsymbol{v})
\end{array}\right],\label{eq:F}
\end{equation}
where we have denoted the tensor product by using the $\otimes$ symbol.
The most common approach to close these equations (see e.g. \citep{reynolds2010})
is to supplement them with the following equation of state
\[
e=\frac{p}{\gamma-1}+\frac{\rho}{2}v^{2}+B^{2}.
\]

In this paper, as in \citep{reynolds2006} and \citep{reynolds2010},
we will consider a slightly more general class of equations which,
in addition to the dynamics discussed so far, includes dissipative
effects (due to particle collisions in the plasma). To that end the
hyperbolic flux vector $F(U)$ is extended in \citep{reynolds2006}
by a diffusive part given by 
\begin{equation}
F_{d}(U)=\left[\begin{array}{c}
0\\
\mu\boldsymbol{\tau}\\
\eta\left(\nabla\boldsymbol{B}-(\nabla\boldsymbol{B})^{\mathrm{T}}\right)\\
\mu\boldsymbol{\tau}\cdot\boldsymbol{v}+\frac{\gamma\mu\kappa}{\gamma-1}\nabla e+\eta\left(\tfrac{1}{2}\nabla(\boldsymbol{B}\cdot\boldsymbol{B})-\boldsymbol{B}\cdot(\nabla\boldsymbol{B})^{\mathrm{T}}\right)
\end{array}\right],\label{eq:Fd-1}
\end{equation}

where
\[
\boldsymbol{\tau}=\nabla\boldsymbol{v}+(\nabla\boldsymbol{v})^{\mathrm{T}}-\tfrac{2}{3}(\nabla\cdot\boldsymbol{\tau})I.
\]
Note that $F_{d}(U)$ as stated here is formulated in terms of the
(dimensionless) viscosity $\mu=Re^{-1}$, the (dimensionless) resistivity
$\eta=S^{-1}$, and the (dimensionless) thermal conductivity $\kappa=Pr^{-1}$,
where $Re$ is the Reynolds number, $S$ the Lundquist number, and
$Pr$ the Prandtl number. In all the simulations conducted, the heat
capacity ratio $\gamma$ is set to $5/3$ which corresponds to a monoatomic
ideal gas.

Thus, the form of the resistive MHD equations used for the spatial
discretization is

\[
\frac{\partial U}{\partial t}+\nabla\cdot F(U)=\nabla\cdot F_{d}(U)
\]
with $F(U)$ given by equation (\ref{eq:F}) and where $F_{d}(U)$
is given by equation (\ref{eq:Fd-1}).

For comparison and validation we have used the Fortran code developed
in \citep{reynolds2006}. This code has been used to conduct plasma
physics simulations (see e.g. \citep{reynolds2006}, \citep{reynolds2010},
and \citep{reynolds2012-2QOYZWMPACZAJ2MABGMOZ6CCPY}) as well as to
construct more efficient preconditioners in the context of implicit
time integrators (see e.g.~\citep{reynolds2010}).

Our implementation (as is the case for D. R. Reynolds' Fortran code)
assumes a finite difference or finite volume method where a single
value is stored in each cell. Thus, in each cell we store the value
of the density $\rho$, the fluid velocity $\boldsymbol{v}$, the
magnetic field $\boldsymbol{B}$, the energy $e$ (but not the pressure
$p$). In principle it is possible to implement any space discretization
satisfying the constraints outlined above. For the numerical simulations
conducted in section \ref{sec:Numerical-examples} we have implemented
the 2.5-dimensional case (that is, the state variables do only depend
on the $x$- and $y$-direction but both the velocity $\boldsymbol{v}(x,y)$
as well as the magnetic field $\boldsymbol{B}(x,y)$ are three-dimensional
vectors) using a classic centered stencil for the divergence; that
is, for each vector $G$, corresponding to the flux of a given (scalar)
state variable, we compute the following approximation of the divergence
\[
\nabla\cdot G(U)\approx\frac{G(U_{i+1,j})-G(U_{i-1,j})}{2h}+\frac{G(U_{i,j+1})-G(U_{i,j-1})}{2h},
\]

where $i$ and $j$ are the cell indices in the $x$- and $y$-directions,
respectively and $h$ is the cell size. To evaluate the spatial derivatives
present in the diffusion vector $F_{d}$ we employ the classic centered
stencil in the interior of the domain and a one-sided stencil of order
$2$ near the boundary. 

In order to integrate the resulting system of ordinary differential
equations in time, we use the CVODE software package. The CVODE library
employs a multistep scheme based on the backward differentiation formulas
(BDF).

\section{Design principles\label{sec:Design-principles}}

In this section we describe, in some detail, some of the design considerations
outlined in the introduction. For the development of the program we
have used the clang C++ compiler as it provides more discernible error
messages. However, in preliminary tests it was observed that the clang
C++ compiler is, at least in some tests, significantly slower than
both the gcc compiler and the Intel C++ compiler. Therefore, and because
of the fact that the clang C++ compiler is usually not present on
most supercomputers, in what follows, we provide only time measurements
for the gcc compiler and the Intel C++ compiler.

\subsection{The choice of an array class}

Compared to Fortran (where usually the single and multi-dimensional
arrays built into the syntax of the language are the best choice)
there are a number of options available in the C++ programming language.
The (traditional) C-style arrays have only limited support for multi-dimensional
arrays and additional information (such as the length of the array)
has to be stored separately. The latter (similar to the way arrays
are handled in Fortran) increases code complexity while the former
(contrary to how arrays are handled in Fortran) requires the user
to develop a custom notation for array indexing.%
\footnote{In the C programming language macros are usually employed to accomplish
this behavior.%
}

In C++ two classes have been included in the standard library that
can be used to replace dynamically allocated arrays (i.e. arrays where
the number of elements is not known at compile time). The \texttt{vector}
class and the \texttt{valarray }class. Both classes handle the allocation
of memory, store the dimension of the array (which can be accessed
by the \texttt{size} method), and allow access to elements by using
the \texttt{{[}{]}} operator. Furthermore, both classes maintain compatibility
to the standard C-style arrays in that all data are stored sequentially
in memory.

The reason the \texttt{valarray} class has been introduced in C++
is to achieve performance comparable to Fortran. The idea is that
since the \texttt{valarray} class will never hold any pointer type
(thus the name value array) it can be implemented more efficiently
than is the case for a generic container (such as the \texttt{vector}
class). More detailed information can be found in \citep{stroustrup1997}.

In Table \ref{fig:array-comparison} the computation of the inviscid
flux for a one-dimensional array of a given size is timed for the
C-style array, the \texttt{vector} class, and the \texttt{valarray}
class. We observe identical results (within measurement errors) with
two exceptions. For the gcc compiler the \texttt{valarray }class is
faster by about 50\% as compared to all other implementations and
for the Intel C++ compiler the Fortran implementation is faster as
compared to all other implementations by approximately the same margin.
In our implementation we have chosen to use the \texttt{vector} class
as there is no clear speed advantage of the \texttt{valarray }class
(in fact in one configuration the \texttt{valarray} class is faster
and in another configuration the \texttt{vector} class is faster)
and the \texttt{vector} class works better with the concept of iterators. 

\begin{table}
\begin{centering}
\begin{tabular}{c|ccll}
 & gcc (LEO III) & icc (LEO III) & gcc (VSC-2) & icc (VSC-2)\tabularnewline
\hline 
C-style & $2.313(93)$ ms & $2.544(99)$ ms & $3.909(71)$ ms & $5.76(25)$ ms\tabularnewline
vector & $2.279(31)$ ms & $2.529(44)$ ms & $3.079(53)$ ms & $6.107(28)$ ms\tabularnewline
valarray & $1.429(52)$ ms & $2.577(45)$ ms & $3.433(59)$ ms & $5.883(43)$ ms\tabularnewline
Fortran & $2.223(27)$ ms & $1.452(45)$ ms & $4.50(21)$ ms & $3.46(12)$ ms\tabularnewline
\end{tabular}
\par\end{centering}

\caption{The inviscid flux for an array of size $n=10^{5}$ is computed. The
measured times are shown (the estimated standard deviations are computed
using $10^{4}$ repetitions and are shown in parenthesis). The simulations
are conducted on the LEO III (Intel Xeon X5650 with gcc 4.4.6, icc
12.1, and Boost 1.51.0) and the VSC-2 supercomputer (AMD Opteron 6132
HE with gcc 4.4.7, icc 14.0.2, Boost 1.55.0). In both cases the \texttt{-O3}
compiler flag is used. \label{fig:array-comparison}}
\end{table}

In the \texttt{C++} standard libraries no class has been included
to facilitate multi-dimensional arrays. However, the Boost library
includes the \texttt{multi\_array} class. It overloads the \texttt{{[}{]}{[}{]}}
and the \texttt{(,)} operator to access elements and provides a mean
to obtain lower-dimensional slices without copying any data (by returning
an iterator). Especially the latter does contribute to code simplicity
as no stride value has to be taken into account. Furthermore, the
\texttt{multi\_array} class maintains memory compatibility with C-style
arrays (the corresponding one-dimensional array can be stored either
in row-major or in column-major format). Of course, performance is
dependent on the order of iteration. An iteration with stride $1$
will, in general, result in better performance compared to an iteration
with stride $N$ (where $N$ is the size of the array in a single
dimension). This has to be taken into account; however, the same is
true for any array implementation. 

In Table \ref{fig:multiarray-comparison} the computation of the inviscid
flux for a two-dimensional array of a given size is timed for the
C-style array and the Boost \texttt{multi\_array} class. We observe
that in this test the performance penalty for using the \texttt{multi\_array}
class as compared to a C-style array is below 10\% if the gcc compiler
is used. Let us further note that the Intel C++ compiler is significantly
slower than the gcc compiler in all of these tests (although the Intel
C++ compiler achieves similar performance compared to the gcc compiler
in case of the Fortran implementation). Let us further note that a
convenient feature of the \texttt{multi\_array} class is that bound
checking can be enabled while the application is being developed (in
the performance measurements conducted here these checks, as is usually
the case in a release build, have been disabled by defining \texttt{BOOST\_DISABLE\_ASSERTS})

\begin{table}
\begin{centering}
Stride $1$\\[1em]
\par\end{centering}

\begin{centering}
\begin{tabular}{c|llll}
 & gcc (LEO III) & icc (LEO III) & gcc (VSC-2) & icc (VSC-2)\tabularnewline
\hline 
C-style & $1.759(32)$ ms & $2.63(10)$ ms & $3.631(69)$ ms & $5.025(22)$ ms\tabularnewline
Fortran & $1.785(82)$ ms & $1.727(73)$ ms & $3.45(33)$ ms & $3.250(26)$ ms\tabularnewline
muti\_array & $1.911(73)$ ms & $3.153(64)$ ms & $3.455(95)$ ms & $6.124(17)$ ms\tabularnewline
\end{tabular}
\par\end{centering}

\begin{centering}
$\phantom{}$\\[0.5em]Stride $N$\\[1em]
\par\end{centering}

\begin{centering}
\begin{tabular}{ccccc}
 & gcc (LEO III) & icc (LEO III) & gcc (VSC-2) & icc (VSC-2)\tabularnewline
\hline 
C-style & $4.25(41)$ ms & $6.96(86)$ ms & $5.899(38)$ ms & $12.632(31)$ ms\tabularnewline
Fortran & $3.76(43)$ ms & $4.37(72)$ ms & $4.1926(87)$ ms & $3.841(15)$ ms\tabularnewline
multi\_array & $4.61(76)$ ms & $7.45(78)$ ms & $6.080(17)$ ms & $13.535(24)$ ms\tabularnewline
\end{tabular}
\par\end{centering}

\caption{The inviscid flux for an array of size $n=N^{2}=316\cdot316\approx10^{5}$
is computed. The measured times are shown (the estimated standard
deviations are computed using $10^{4}$ repetitions and are shown
in parenthesis). The simulations are conducted on the LEO III (Intel
Xeon X5650 with gcc 4.4.6, icc 12.1, and Boost 1.51.0) and the VSC-2
supercomputer (AMD Opteron 6132 HE with gcc 4.4.7, icc 14.0.2, Boost
1.55.0). In both cases the \texttt{-O3} compiler flag is used. \label{fig:multiarray-comparison}}
\end{table}

\subsection{Mapping state variables to physically relevant variables\label{sub:Mapping-state-variables}}

In numerical simulations the phase space is usually represented as
an array of double precision floating point numbers. However, the
physical variables associated to this quantities are usually referred
to by a well known designation (such as $\rho$ for the density or
$\boldsymbol{B}$ for the magnetic field). From this consideration
two issues arise. First, it is not desirable to access physical quantities
purely by their array index. For example, code fragments similar to
\begin{lstlisting}[language=Matlab]
y(i,1)=x(i,1)*x(i,3)
\end{lstlisting}
are used in many Fortran, Matlab, and C implementations to compute
the mass flux with respect to the $y$-direction. Obviously code written
in this fashion is not easily discernible even if it is known that
the first component represents the density and that the third component
represents the velocity in the $y$-direction. Second, some quantities
can be considered as vector types (for example the magnetic field).
Thus, some way to iterate over such variables has to be provided.
A common way to implement this in C++ (or C) is to define an appropriate
\texttt{struct.} In our case the state space would be represented
by

\begin{lstlisting}[language={C++}]
struct state {
  double rho;
  array<double,3> v;
  array<double,3> B;
  double e;
};
\end{lstlisting}
The issue with this implementation is that to implement collective
operation on the entire state (such as component-wise addition and
scalar multiplication) is either somewhat tedious as all the variables
have been summed separately or does rely on techniques which are not
type-safe (such as casting the struct to a double pointer). 

Therefore, we have chosen to implement the state as a simple $8$-dimensional
array and to supply appropriate functions to access the physical variables.

\begin{lstlisting}[language={C++}]
struct state {
	array<double, 8> data;
};

double& rho(state& s) {
  return s.data[0]; 
}
double& v(state& s, size_t dim) {
  return s.data[1+dim]; 
}
...
\end{lstlisting}
As can be seen from the code snippet, we employ a feature of C++ that
allows us to return a value as a reference from a function. This in
fact allows us to set as well as access (get) the desired values.
The computation of the flux from the beginning of the section in this
notation reads as

\begin{lstlisting}[language={C++}]
rho(y)=rho(x)*v(x,1);
\end{lstlisting}
We expect the compiler to inline this calls and thus no performance
overhead should be incurred. We have conducted a number of numerical
simulations which confirms this expectation. 

The inviscid flux given in equation (\ref{eq:F}) can be implemented
as follows

\begin{lstlisting}[language={C++}]
state flux_inviscid(size_t dim, state s) {
  state r;
  rho(r) = -rho(s)*v(s, dim);
  double Bv = 0.0;
  for (int i = 0; i < 3; i++) {
    v(r, i) = -rho(s)*v(s,i)*v(s, dim) 
      -(pressure(s) + 0.5*Bsq(s))*(i==dim) + B(s, i)*B(s, dim);
    B(r, i) = v(s, i)*B(s, dim) - B(s, i)*v(s, dim);
    Bv += B(s, i)*v(s, i);
  }
  e(r) = -(e(s)+pressure(s)+0.5*Bsq(s))*v(s, dim) 
       + Bv*B(s, dim);
  return r;
}
\end{lstlisting}
This form greatly increases readability and reflects the mathematical
formulation of the MHD equation (the equation of state is implemented
in the \texttt{pressure} function).

\subsection{Boundary conditions}

One aspect that frequently introduces mistakes in numerical code is
the handling of boundary conditions. There are two commons ways to
implement boundary conditions. First, they can be incorporated directly
into the stencil. This approach, however, has the disadvantage that
the stencil needs to be modified every time different boundary conditions
are considered. In addition, a differentiation between different parts
of the boundary is often necessary. The second approach employs so-called
ghost cells, which are appended to the regular domain. This alleviates
the issues discussed above but requires the extension of the domain.
This in turn makes it often necessary to extract the physical domain
in other parts of the code (for example in a parallel setting where
the values in the interior of the domain have to be passed to a linear
or nonlinear solver). Let us also note that these considerations are
not only important at the physical boundary but also at the boundaries
that exist between different nodes in a MPI parallelization, for example.

We have chosen to store the (inner) domain separately from the ghost
cells (which form a halo region). An iterator is then defined which
allows the developer implementing a stencil to only consider it in
the interior of the domain. The boundary conditions are then automatically
satisfied if the ghost cells are set appropriately. Also an implementation
of additional boundary conditions would not require any modification
of the stencil code. In Table \ref{tab:boundary} we have considered
the impact on performance of this implementation for the classic centered
difference stencil. For the gcc compiler no performance penalty can
be observed while a small decrease in performance is observed for
the Intel C++ compiler. Let us further note that the gcc compiler
seems to perform consistently better than the Intel C++ compiler in
this test.

\begin{table}
\begin{centering}
\begin{tabular}{c|llll}
 & gcc (LEO III) & icc (LEO III) & gcc (VSC-2) & icc (VSC-2)\tabularnewline
\hline 
Direct impl. & $1.086(80)$ ms & $1.282(73)$ ms & $2.49(14)$ ms & $2.47(13)$ ms\tabularnewline
Boundary iter. & $1.085(48)$ ms & $1.411(54)$ ms & $2.234(73)$ ms & $2.637(23)$ ms\tabularnewline
\end{tabular}
\par\end{centering}

\caption{The first derivative for an array of size $n=10^{5}$ is computed
by using the classic centered stencil. The measured times are shown
(the estimated standard deviations are computed using $10^{4}$ repetitions
and are shown in parenthesis). The simulations are conducted on the
LEO III (Intel Xeon X5650 with gcc 4.4.6, icc 12.1, and Boost 1.51.0)
and the VSC-2 supercomputer (AMD Opteron 6132 HE with gcc 4.4.7, icc
14.0.2, Boost 1.55.0). In both cases the \texttt{-O3} compiler flag
is used. \label{tab:boundary}}
\end{table}

\subsection{std::transform\label{sub:functional}}

A common practice in scientific codes is to apply a function to each
element in an array. In this instance one can use either one or multiple
for loops to iterate over all the elements or use the \texttt{transform}
function included in the \texttt{C++} standard library. The latter
conceivably has two advantages. First, it frees us from explicitly
taken care of the array indices and we are thus able to write code
that works independently of the dimension of the array under consideration.
Second, the compiler can take advantage of the information that the
operations for each element of the array are independent and possibly
generate more efficient code. 

In the numerical computations we have conducted (the result of which
are shown in Table \ref{tab:transform}) we observe no clear advantage
for using a for loop or the transform function. In fact, depending
on the configuration either method can yield superior performance).
As expected the \texttt{multi\_array} class shows almost the same
performance as the one-dimensional containers (as there is no difference
in the manner the data are stored in memory). It is further interesting
to note that the gcc compiler outperforms the Intel C++ compiler in
all of these tests.

\begin{table}
\begin{centering}
\begin{tabular}{l|llll}
 & gcc (LEO III) & icc (LEO III) & gcc (VSC-2) & icc (VSC-2)\tabularnewline
\hline 
vector (F) & $2.312(18)$ ms & $2.559(24)$ ms & $3.184(57)$ ms & $6.107(28)$ ms\tabularnewline
vector (T) & $2.273(20)$ ms & $2.515(22)$ ms & $4.51(15)$ ms & $5.645(27)$ ms\tabularnewline
multi\_array (T) & $2.346(67)$ ms & $2.509(29)$ ms & $4.11(11)$ ms & $5.789(12)$ ms\tabularnewline
\end{tabular}
\par\end{centering}

\caption{The inviscid flux for an array of size $n=10^{5}$ is computed using
the vector class (where F denotes that we are using a for loop and
T denotes that we employ the transform function found in the C++ standard
library). The measured times are shown (the estimated standard deviations
are computed using $10^{4}$ repetitions and are shown in parenthesis).
The simulations are conducted on the LEO III (Intel Xeon X5650 with
gcc 4.4.6, icc 12.1, and Boost 1.51.0) and the VSC-2 supercomputer
(AMD Opteron 6132 HE with gcc 4.4.7, icc 14.0.2, Boost 1.55.0). In
both cases the \texttt{-O3} compiler flag is used. \label{tab:transform}}
\end{table}

\section{Numerical examples\label{sec:Numerical-examples}}

\begin{figure}
\begin{centering}
\includegraphics[width=12cm]{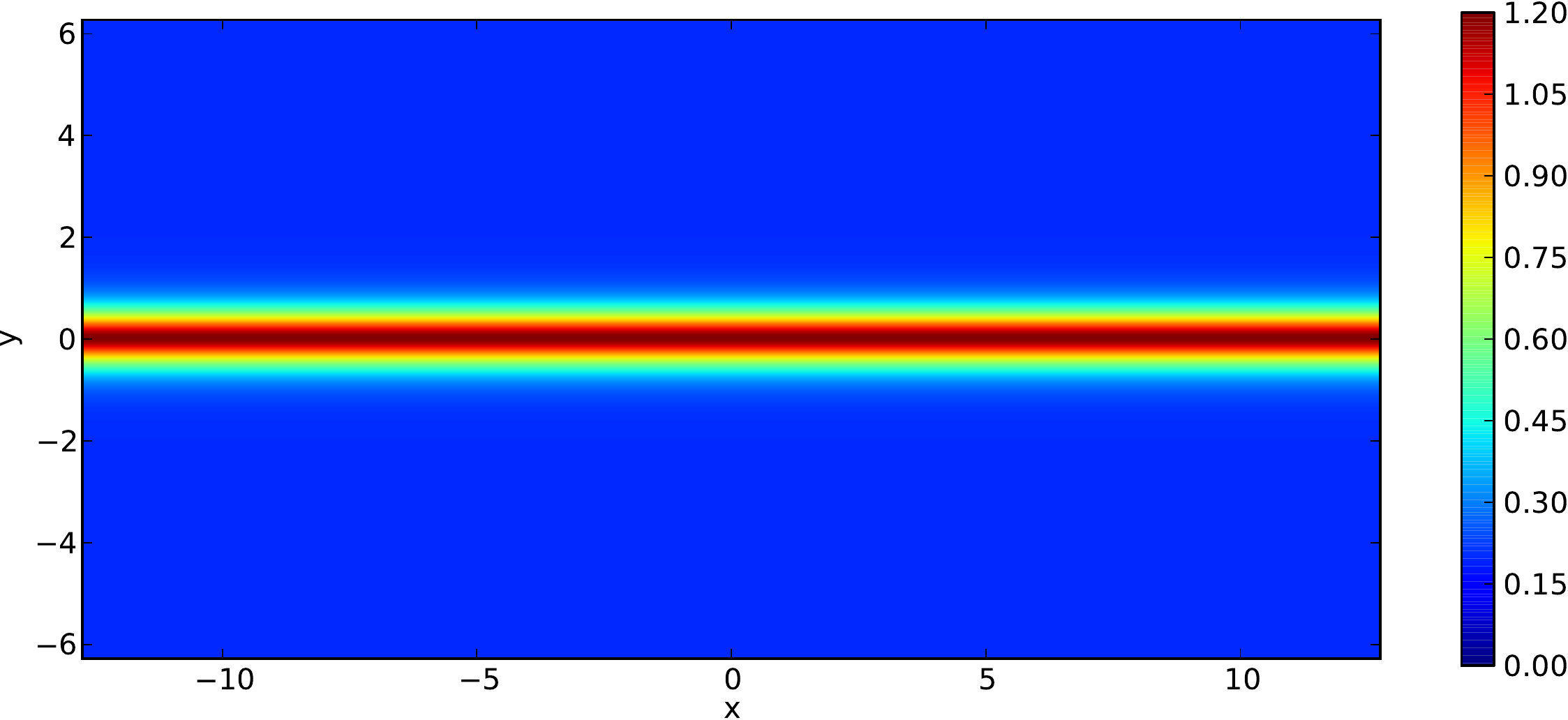}
\par\end{centering}

\begin{centering}
\includegraphics[width=12cm]{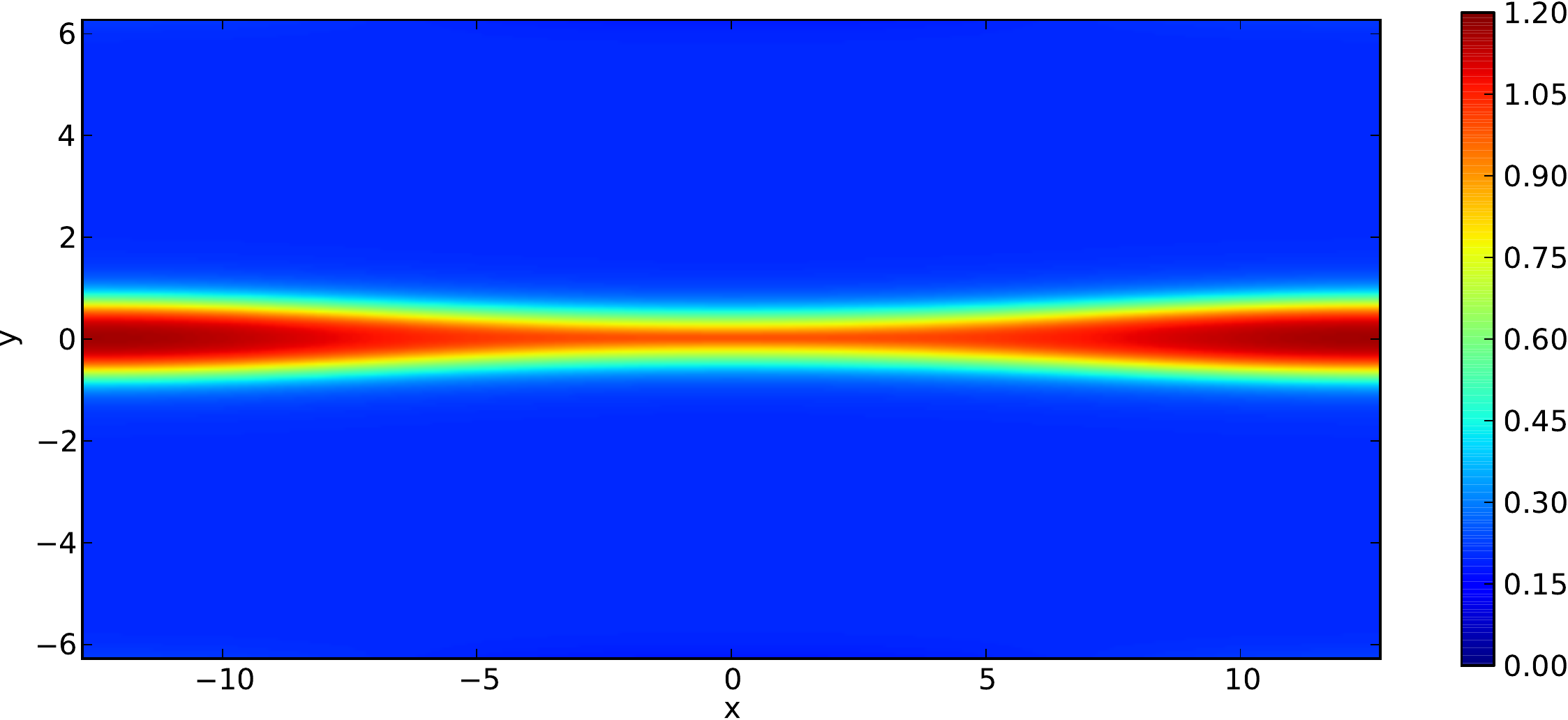}
\par\end{centering}

\begin{centering}
\includegraphics[width=12cm]{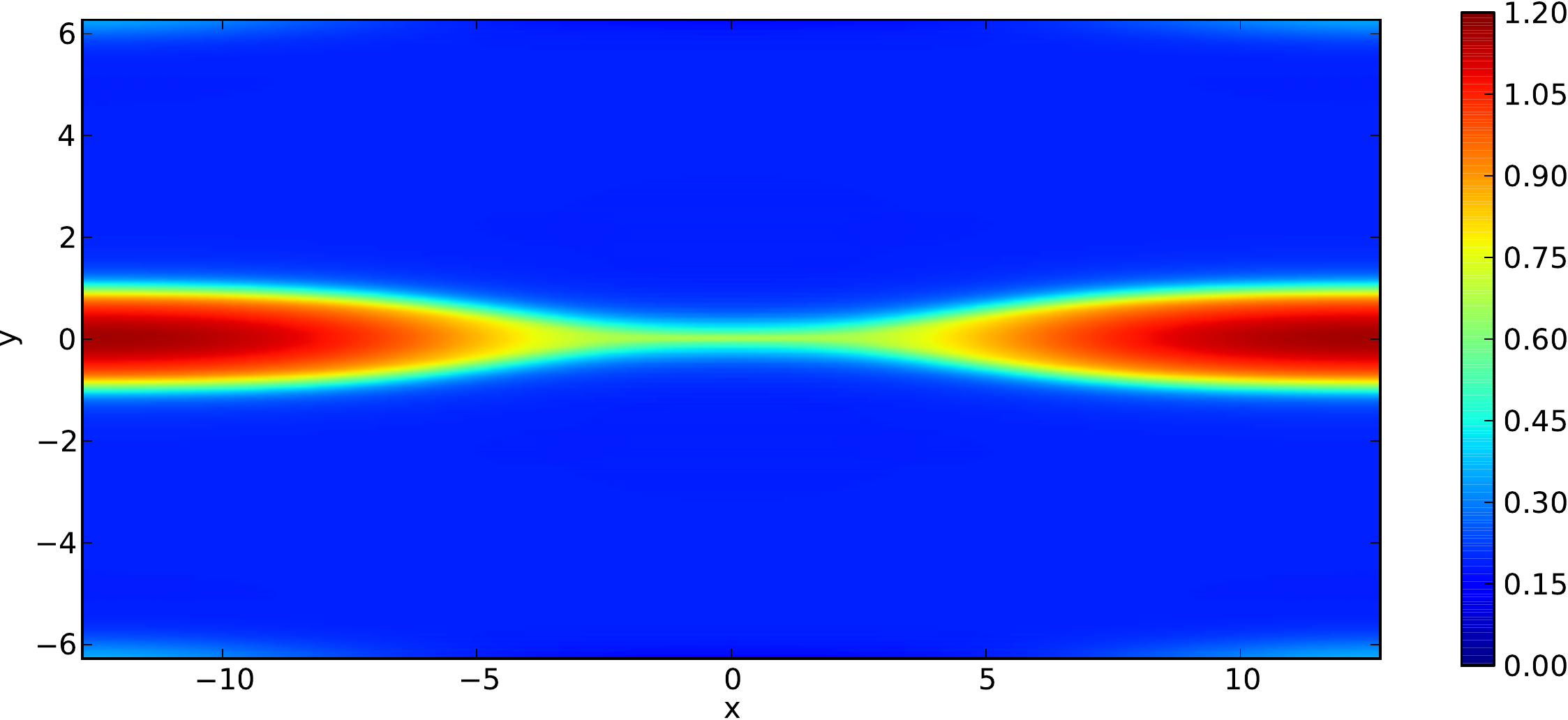}
\par\end{centering}

\caption{The density $\rho$ is shown at times $t=0$ (top), $t=75$ (center),
and $t=150$ (bottom). In both space directions $256$ grid points
are employed. The dimensionless parameters are chosen as follows:
$\mu=10^{-3}$, $\eta=10^{-3}$, and $\kappa=10^{-2}$. \label{fig:recon}}
\end{figure}

To test and validate the code developed, we have run a number of simulations
and compared them to the results produced by D. R. Reynolds' Fortran
code (as described in \citep{reynolds2006} and \citep{reynolds2010}).
The first example we consider here is the so-called reconnection problem.
This example models the reconnection of magnetic field lines. That
is, as an initial value we impose an almost discontinuous magnetic
field which points in opposite directions on different sides of the
discontinuity. The dissipation included in the resistive MHD model
then causes a reconnection of the field lines. The numerical results
(shown in Figure \ref{fig:recon}) display the characteristic broadening/thinning
of the density (see, for example, \citep{reynolds2006}).

\begin{figure}
\begin{centering}
\includegraphics[width=12cm]{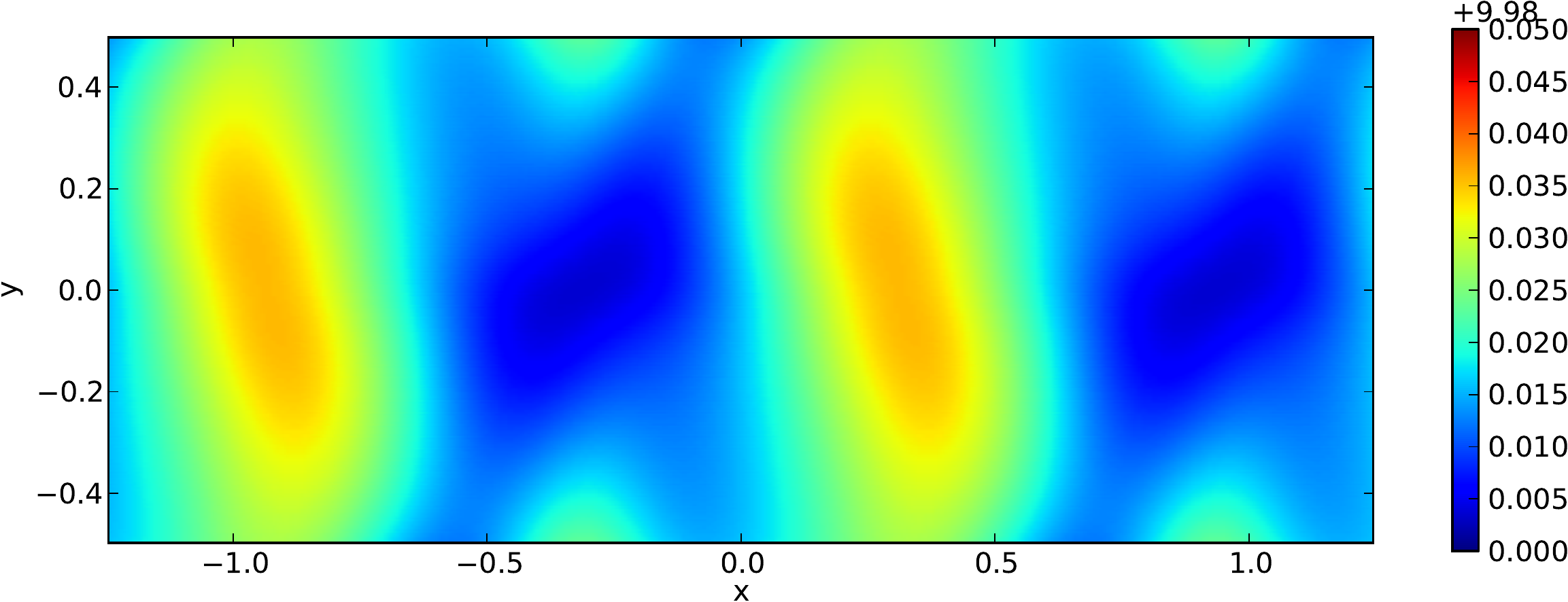}
\par\end{centering}

\begin{centering}
\includegraphics[width=12cm]{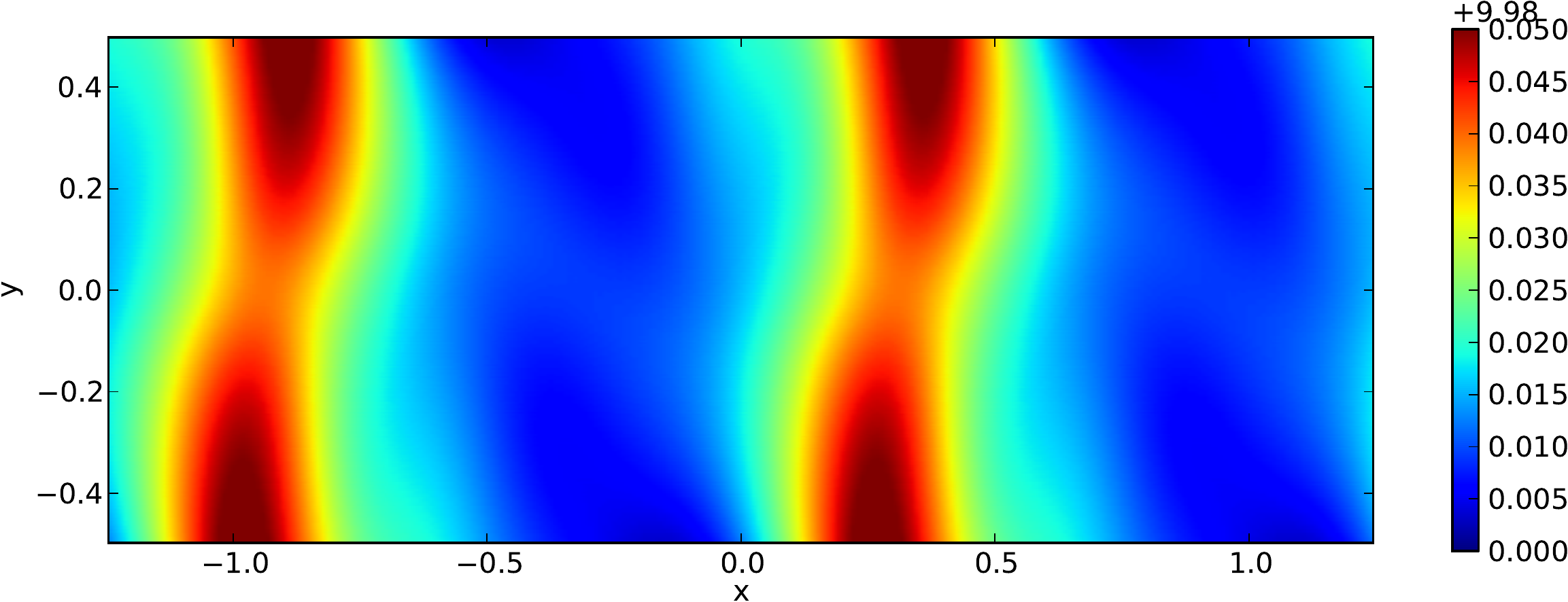}
\par\end{centering}

\caption{The magnitude of the magnetic field in the $z$-direction at times
$t=1$, and $t=2$ is shown. The initial value (i.e. the value at
time $t=0$) consists of a homogeneous field equal in strength to
$10$. In both space directions $256$ grid points are employed. The
dimensionless parameters are chosen as follows: $\mu=10^{-3}$, $\eta=10^{-3}$,
and $\kappa=10^{-2}$. \label{fig:KH}}
\end{figure}

The second example we consider is the so-called Kelvin--Helmholtz
instability. That is, we consider a shear flow in a constant magnetic
field in which the velocity field is perturbed. This leads to a time
evolution with highly inhomogeneous magnetic fields. The results of
the numerical simulation are shown in Figure \ref{fig:KH}.

We have chosen to use Boost.MPI to parallelize our code to multiple
nodes of the cluster under consideration. Boost.MPI provides a convenient
C++ interface to the Message Passing Interface. The only modification
to the existing code that is required to use Boost.MPI is to add a
serialization statement to the \texttt{state} class. This is necessary
as Boost.MPI uses the serialization interface in order to convert
custom data types to a bytestream that is send to different nodes.

In Table \ref{tab:Weak-scaling} a weak scaling test (that is, the
size of the problem is proportional to the number of cores used) for
the reconnection problem is performed. We observe that communication
time necessary to evaluate the right-hand side is not a limiting factor
up to at least 256 cores as most of the increase in the runtime is
due to the increase in computational effort required by the iterative
solver. To reduce this growth in computational effort (to some extent)
a number of advanced preconditioning techniques have been developed
(see, for example, \citep{reynolds2010}). 

\begin{table}
\begin{centering}
\begin{tabular}{c|ccccc}
Cores & Time & RHS eval. & Nonlin. iter. & Linear iter. & Communication time\tabularnewline
\hline 
1 & $0.8$ s & $46$ & 16 & 26 & $0.16$ s\tabularnewline
4 & $2.5$ s & 90 & 34 & 52 & $0.75$ s\tabularnewline
16 & $10.5$ s & 214 & 40 & 170 & $2.4$ s\tabularnewline
64 & $115$ s & 474 & 40 & 430 & $15.9$ s\tabularnewline
256 & $495$ s & 1317 & 44 & 1269 & $20.5$ s\tabularnewline
\end{tabular}
\par\end{centering}

\caption{Weak scaling test for the reconnection problem on the LEO3 HPC system
at the University of Innsbruck (\protect\href{http://www.uibk.ac.at/zid/systeme/hpc-systeme/leo3/}{http://www.uibk.ac.at/zid/systeme/hpc-systeme/leo3/}).
The reconnection problem is integrated up to time $t=1$. The dimensionless
parameters are chosen as follows: $\mu=10^{-3}$, $\eta=10^{-3}$,
and $\kappa=10^{-2}$. To each core we assign a slice of the domain
with $128$ grid points in both directions. The measured communication
time includes both the time necessary to prepare the data as well
as the transmission time. \label{tab:Weak-scaling}}
\end{table}

In Table \ref{tab:Weak-scaling-KH} a weak scaling test for the Kelvin--Helmholtz
instability is performed. We draw the same conclusions as in the case
of the reconnection problem.

\begin{table}
\begin{centering}
\begin{tabular}{c|ccccc}
Cores & Time & RHS eval. & Nonlin. iter. & Linear iter. & Communication time\tabularnewline
\hline 
1 & $7.5$ s & 419 & 202 & 213 & $1.6$ s\tabularnewline
4 & $14$ s & 506 & 167 & 335 & $4$ s\tabularnewline
16 & $39$ s & 634 & 80 & 550 & $12$ s\tabularnewline
64 & $363$ s & 1479 & 104 & 1371 & $55$s\tabularnewline
256 & $1273$ s & 3431 & 112 & 3315 & $55$ s\tabularnewline
\end{tabular}
\par\end{centering}

\caption{Weak scaling test for the reconnection problem on the LEO3 HPC system
at the University of Innsbruck (\protect\href{http://www.uibk.ac.at/zid/systeme/hpc-systeme/leo3/}{http://www.uibk.ac.at/zid/systeme/hpc-systeme/leo3/}).
The Kelvin--Helmholtz instability is integrated up to time $t=0.05$.
To each core we assign a slice of the domain with $128$ grid points
in both directions. The measured communication time includes both
the time necessary to prepare the data as well as the transmission
time. \label{tab:Weak-scaling-KH}}
\end{table}

However, what is more relevant with respect to the implementation
considered here is how the communication cost increases as we increase
the number of cores. In the weak scaling tests conducted above this
behavior can not be easily observed as the total computational cost
also increases as the problem size is increased. However, it is difficult
to conduct a strong scaling test using the implicit time integrator
as the memory available per core (on both HPC systems under consideration
here) does not allow us to run a problem of significant size on a
single node. Therefore, in Table \ref{tab:strong-scaling} we show
the results of a strong scaling test using the explicit Euler method
as the time integrator. 

\begin{table}
\begin{centering}
\begin{tabular}{c|ccc}
Cores & Time (1024) & Speedup (1024) & Communication time\tabularnewline
\hline 
1 & $126$ s & - & -\tabularnewline
16 & $13$ s & $10$ & $3.2$ s\tabularnewline
64 & $10$ s & $13$ & $8.1$ s\tabularnewline
\end{tabular}
\par\end{centering}

\caption{Strong scaling test for the reconnection problem on the LEO3 HPC system
at the University of Innsbruck (\protect\href{http://www.uibk.ac.at/zid/systeme/hpc-systeme/leo3/}{http://www.uibk.ac.at/zid/systeme/hpc-systeme/leo3/}).
The reconnection problem is integrated up to time $t=0.1$. The dimensionless
parameters are chosen as follows: $\mu=10^{-3}$, $\eta=10^{-3}$,
and $\kappa=10^{-2}$. The domain consists of $1024$ grid points
in both directions. The measured communication time includes both
the time necessary to prepare the data as well as the transmission
time. \label{tab:strong-scaling}}
\end{table}

Furthermore, we have compared the time it takes to compute the flux
vector for D. R. Reynolds' Fortran implementation and our code. For
the Fortran implementation we measure approximately $12$ ms (using
the gcc Fortran compiler) and $10$ ms (using the Intel Fortran compiler).
For our C++ implementation we measure approximately $7.5$ ms (using
the gcc compiler) and $8.5$ ms (using the Intel C++ compiler). These
simulations have been conducted on the LEO3 HPC system and are consistent
with the performance measurements presented in section \ref{sec:Design-principles}.
The results show a performance advantage of our C++ code of approximately
30\%. Note, however, that caution is warranted in comparing these
two implementations as there are small differences with respect to
both the implementation and the numerical method that is implemented.
However, what the comparison does show is that the performance achieved
by the C++ code is certainly comparable to that of the Fortran code.
Certainly there are a number of optimizations that, if carried out,
could improve the performance of both codes. We should also emphasize
at this point that the C++ implementation is certainly significantly
shorter (measured in number of lines of code) compared to the Fortran
implementation.

\section{Documentation\label{sec:Documentation}}

The purpose of this section is to provide an overview of the C++ code
discussed in this paper. All the routines which are not directly related
to the MHD problem or the numerical scheme have been implemented in
C++ header files (with the ending .hpp). These header files, as is
necessary for template based programming, contain both the class and
function declarations as well as their implementation. The files are
\begin{description}
\item [{state.hpp}] This file includes the \texttt{state} class and all
the helper functions that are used to manipulate it. This includes
the functions discussed in section~\ref{sub:Mapping-state-variables}.
\item [{boundary.hpp}] This file includes the \texttt{boundary2d} class
which stores the boundary conditions (both left/right and top/bottom).
The \texttt{boundary2d} class includes methods to obtain iterators
over the boundary elements.
\item [{domain.hpp}] This file includes two classes. The \texttt{interval}
class which is simply used to store the size of a given domain (in
a single direction) and the \texttt{domain} class. The \texttt{domain}
class stores a two-dimensional array with the values corresponding
to the part of the domain residing on the corresponding processor.
In addition, it includes a method to construct an initial approximation
(from a function depending on both the $x$- and $y$-direction),
a method to obtain a one-dimensional slice of the data, and a function
to write the data to disk.
\item [{stencil\_array.hpp}] This file includes the \texttt{stencil\_array}
class which constructs a (virtual) array from an iterator to a slice
of the interior domain and the appropriate boundary elements. The
\texttt{{[}{]}} operator is overloaded in order to provide access
to these (virtual) one-dimensional array. 
\item [{stencil.hpp}] This file includes the single function \texttt{apply\_stencil}
which applies a stencil (in the $x$-direction if the parameter dim
is set to $0$ or the $y$- direction if the parameter dim is set
to $1$) using the domain data \texttt{dom} and the boundary data
\texttt{bdr}. The output is stored in \texttt{out}. The template parameter
\texttt{stencil} (of type FUNC) is a function which computes the desired
stencil giving a reference to an input and an output stencil array
(these can be accessed as generic one-dimensional arrays as discussed
above).
\item [{mpi.hpp}] This file includes the \texttt{parallel\_mpi} class which
is responsible for mapping the MPI rank to the two-dimensional indices
in the grid used for conducting the computation. Furthermore, it provides
helper methods to check which physical boundaries the current core
needs to handle.
\item [{cvode.hpp}] This file includes two classes: \texttt{wrapper} and
\texttt{cvode\_wrapper}. The former is a base class with a single
virtual method. The second implements this method in order to conduct
a time step using the CVODE library. All the code that interfaces
with the CVODE library is included in this file.
\item [{timer.hpp}] This files includes the \texttt{timer} class which
is used in order to time the code.
\end{description}
All of these files combined contain about 500 lines of source code.
Their primary purpose is to separate the common tasks from the specific
implementation of a given set of equations, the numerical scheme used,
and the communication required. These remaining functions are all
implemented in the mhd.cpp file which contains approximately 300 lines
of code (of which about 100 lines of code are used to parse the command
line arguments; the Boost.Program\_options library is used for that
purpose).

The program accepts a number of command line option. Using the command

\begin{lstlisting}[language=bash]
./mhd --help
\end{lstlisting}
a list including a more detailed description can be obtained. For
the two examples introduced in the previous section the commands

\begin{lstlisting}[language=bash]
mpirun -np 4 ./mhd -m cvode -p recon -T 150 \
       --timestep 25 --nx 256 --ny 256
\end{lstlisting}
and

\begin{lstlisting}[language=bash]
mpirun -np 4 ./mhd -m cvode -p kh -T 2 \
       --timestep 0.1 --nx 256 --ny 256
\end{lstlisting}
have been used. By default the workload is equally distributed among
the two directions.

\section{Conclusion \& Outlook\label{sec:Conclusion-&-Outlook}}

We have outlined the design rationals upon which the implementation
of our resistive magnetohydrodynamics solver is based. From our performance
measurements it should be clear that this greatly improves code readability
while it (at most) incurs a minimal performance penalty. In fact,
the implementation compares very favorable to an existing Fortran
implementation from the literature. Beyond the \texttt{C++} standard
library we only use the Boost library as well as CVODE (in order to
conduct the time integration). It is quite surprising to us that (except
for the Fortran benchmarks) the Intel C compiler was outperformed
by the gcc compiler in almost every test. This behavior is consistent
across at least two HPC systems. 

The framework developed here can conceivably be used for problems
different from MHD and additional time integration libraries can be
included (an implementation that includes the, as of yet unpublished,
EPIC exponential integrator library is currently in work).

\section{Acknowledgement}

We would like to take the opportunity to thank D. R. Reynolds for
providing the code of his Fortran MHD solver. This has greatly aided
us in the testing and validation process of the resistive MHD solver
described here.

\bibliographystyle{plainnat}
\bibliography{papers}

\begin{thebibliography}{8}
\providecommand{\natexlab}[1]{#1}
\providecommand{\url}[1]{\texttt{#1}}
\expandafter\ifx\csname urlstyle\endcsname\relax
  \providecommand{\doi}[1]{doi: #1}\else
  \providecommand{\doi}{doi: \begingroup \urlstyle{rm}\Url}\fi

\bibitem[boo()]{boost}
{BOOST C++ Libraries}.
\newblock URL \url{http://www.boost.org}.
\newblock Last retrieved \today.

\bibitem[Cary et~al.(1997)Cary, Shasharina, Cummings, Reynders, and
  Hinker]{cary1997}
J.R. Cary, S.G. Shasharina, J.C. Cummings, J.V.W. Reynders, and P.J. Hinker.
\newblock {Comparison of C++ and Fortran 90 for object-oriented scientific
  programming}.
\newblock \emph{Comput. Phys. Commun.}, 105\penalty0 (1):\penalty0 20--36,
  1997.

\bibitem[Pesch et~al.(2007)Pesch, Bell, Sollie, Ambati, Bokhove, and Van
  Der~Vegt]{pesch2007}
L.~Pesch, A.~Bell, H.~Sollie, V.R. Ambati, O.~Bokhove, and J.J. Van Der~Vegt.
\newblock {hpGEM---A software framework for discontinuous Galerkin finite
  element methods}.
\newblock \emph{ACM T. Math. Software}, 33\penalty0 (4):\penalty0 23, 2007.

\bibitem[Reynolds et~al.(2006)Reynolds, Samtaney, and Woodward]{reynolds2006}
D.R. Reynolds, R.~Samtaney, and C.S. Woodward.
\newblock {A fully implicit numerical method for single-fluid resistive
  magnetohydrodynamics}.
\newblock \emph{J. Comput. Phys.}, 219\penalty0 (1), 2006.

\bibitem[Reynolds et~al.(2010)Reynolds, Samtaney, and Woodward]{reynolds2010}
D.R. Reynolds, R.~Samtaney, and C.S. Woodward.
\newblock {Operator-based preconditioning of stiff hyperbolic systems}.
\newblock \emph{SIAM J. Sci. Comput.}, 32\penalty0 (1), 2010.

\bibitem[Reynolds et~al.(2012)Reynolds, Samtaney, and
  Tiedeman]{reynolds2012-2QOYZWMPACZAJ2MABGMOZ6CCPY}
D.R. Reynolds, R.~Samtaney, and H.C. Tiedeman.
\newblock {A fully implicit Newton--Krylov--Schwarz method for tokamak
  magnetohydrodynamics: Jacobian construction and preconditioner formulation}.
\newblock \emph{Computational Science \& Discovery}, 5\penalty0 (1), 2012.

\bibitem[Stroustrup(1997)]{stroustrup1997}
B.~Stroustrup.
\newblock \emph{{The C++ Programming Language}}.
\newblock Addison Wesley, 3rd edition, 1997.

\bibitem[Veldhuizen(1998)]{veldhuizen1998}
T.L. Veldhuizen.
\newblock {Arrays in Blitz++}.
\newblock In \emph{{Computing in Object-Oriented Parallel Environments, Lecture
  Notes in Computer Science}}, volume 1505, pages 223--230. Springer, 1998.

\end{thebibliography}

\end{document}